\documentstyle[12pt,amssymb]{article}
\input psfig.sty

\newtheorem{theorem}{Theorem} 
\newtheorem{remark}{Remark}

\newenvironment{proof}{{\bf Proof: }}{} 
\newcommand{\R}{{\Bbb R}}
\newcommand{\st}{{\bigm|}}
\newcommand{\opname}[1]{\mathop{\fam0#1}}
\newcommand{\area}{\opname{area}}
\newcommand{\qed}{\relax{\ifhmode\unskip\nobreak\hfill$\Box$\fi
\ifmmode\ifinner\else\hskip5pt\fi\hfill\Box\fi}\relax}
\newcommand{\bq}{\begin{equation}}
\newcommand{\eq}{\end{equation}}

\begin{document}
\title{Average kissing numbers \\ for non-congruent sphere packings}
\author{Greg Kuperberg
\thanks{Supported by an NSF Postdoctoral Fellowship, grant \#DMS-9107908}\\
Oded Schramm
\thanks{Incumbent of the William Z.\ and Eda Bess Novick Career 
Development Chair.  Supported by NSF grant \#DMS-9112150}}
\date{May 13, 1994}
\maketitle

\section{Introduction}

Let $P$ be a packing of $n$ (round) balls in $\R^3$.  (A packing of round
balls, also known as a sphere packing, is a collection of round balls with
disjoint interiors.)  The balls may have different radii.  The average
kissing number of $P$ is defined as $k(P) = 2m/n$, where $m$ is the number of
tangencies between balls in the packing.  Let 
$$k = \sup\{k(P)\st P\mbox{ is a finite packing of balls in }\R^3\}.$$
\begin{theorem} 
$$ 12.566 \approx 666/53 \le k < 8+4\sqrt 3 \approx 14.928.$$
\label{kest}
\end{theorem}

(The appearance of the number of the beast in the lower bound is
purely coincidental.)

The supremal average kissing number $k$ is defined in any dimension, as are 
$k_c$, the supremal average kissing number for congruent ball packing, and
$k_s$, the maximal kissing number for a single ball surrounded by congruent
balls with disjoint interiors.  (Clearly, $k_c \le k$ and $k_c \le k_s$.) It
is interesting that $k$ is always finite, because a large ball can be
surrounded by many small balls in a non-congruent ball packing.  
Nevertheless, a simple argument presented below shows that $k \le 2k_s$ in
every dimension, and clearly $k_s$ is always finite. In two dimensions, an
Euler characteristic argument shows that $k \le 6$, but it is also well-known
that $k_s = k_c = 6$. One might therefore conjecture that $k = k_c$ always,
or at least in dimensions such as 2, 3, 8, and 24 (and conjecturally several
others) in which $k_s = k_c$ \cite{C+S:book}.   Surprisingly, in three
dimensions,  $k > 12$ even though $k_s = k_c = 12$.

\begin{remark} \em No packing $P$ achieves the supremum $k = k(P)$,
because if $P'$ is a translate of $P$ that meets $P$ in only one
point, then $k(P \cup P') > k(P)$.
\label{ksup}
\end{remark}

Let $P=(P_v,v\in V)$ be a packing, where $V$ is some indexing set. The {\it
nerve\/} of $P$ is a combinatorial object that encodes the combinatorics of
the packing.  It is the (abstract) graph $G=(V,E)$ on $V$, where an edge
$\{u,w\}$ appears in $E$ precisely when  $P_u$ and $P_w$ intersect. If $P$ is
a packing of round disks in the plane, then it is easy to see that $G$ is a
planar graph.  Conversely, the circle packing theorem \cite{Koebe:circpac},
states that every finite planar graph is the nerve of some disk packing in
the plane.  This non-trivial theorem has received much attention lately,
mostly because of its surprising relation with complex analysis. (Compare
references \cite{Thurston:book}, \cite{R+S}, and \cite{CdV:conv}.)

Since the nerves of planar disk packings are understood, it is natural to ask
for a description of all graphs that are nerves of ball packings in $\R^3$.
In lieu of a complete characterization, which is probably intractable,
Theorem~\ref{kest} gives a necessary condition on such graphs:
$2|E| < (8+4\sqrt 3)|V|$.


We wish to thank Gil Kalai for a discussion which led to the question of
estimating $k$.

\section{The upper bound}

\begin{theorem}
If $P$ is a finite ball packing in $\R^3$, then $k(P)<8 + 4\sqrt 3$.
\label{ub}
\end{theorem}

As a warm-up, we will show that $k(P)\le 24$.  Let $E$ be the set of unordered
pairs of balls in $P$ that kiss.  Let $r(B)$ be the radius of a ball $B \in
P$.  By a famous result \cite{S+W}, \cite{Leech:13spheres}, it is impossible
for more than 12 unit balls with disjoint interiors to kiss a unit ball $B$.
If $C$ kisses $B$ and $r(C) > 1=r(B)$, then $C$ contains a (unique) unit ball
that kisses $B$.  Thus, in a packing, $B$ cannot kiss more than 12 balls at
least as large as $B$.  Consider a function $f:E \to P$ that assigns to
$\{B,C\}\in E$ the smaller of the balls $B$ and $C$, or either if they are
the same size.  Since $f$ is at most 12 to 1, $|E| \le 12 |P|$. Consequently,
$k(P) = 2|E|/|P| \le 24$.

The proof of Theorem~\ref{ub} is a refinement of this argument.

\begin{proof}
In addition to the above notation, we let $E(B)$ denote the set of $C\in P$
such that $\{B,C\}\in E$.

Let $\rho>1$ be a constant to be determined below. For each ball $B
\in P$, let $S(B)$ be the concentric spherical shell with radius $\rho r(B)$.
For each $B,C \in P$, define
\bq a(B,C) =\frac{\area(C\cap S(B))}{\area(S(B))}. \label{humpty}\eq
Since the interiors of the balls in $P$ are disjoint,  for any $B$,
\bq 1 \ge \sum_{C \in P} a(B,C) \ge \sum_{C\in E(B)} a(B,C).
\label{start}\eq
Summing over $B$,
\bq |P| 
\ge \sum_{\{B,C\} \in E} \left(a(B,C) + a(C,B)\right). \label{d1}\eq

\begin{figure}[ht]
\centerline{\mbox{\psfig{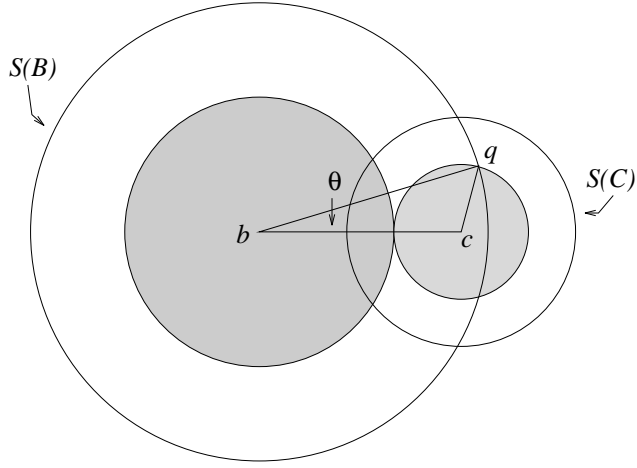}}}
\caption{The intersection of $B$ and $C$ with a plane passing through their
centers.}
\label{bc}
\end{figure}

We will obtain a lower bound on $a(B,C)+a(C,B)$ for two kissing balls $B$
and $C$.  Suppose that $B$ insersects $S(C)$ and $C$ intersects $S(B)$, as
shown in Figure~\ref{bc}.  Let $b$ and $c$ be the centers of $B$ and
$C$.  Let $q$ be a point on the relative boundary in $S(B)$ of the spherical
disk $C \cap S(B)$.  Clearly,
\begin{eqnarray*}
d(b,c) &=& r(B)+r(C)\\
d(b,q) &=& \rho r(B)\\
d(c,q) &=& r(C),
\end{eqnarray*}
where $d(x,y)$ is the distance from $x$ to $y$.  Let $\theta = \angle cbq$ be
the angular radius of $C \cap S(B)$.  By the law of cosines,
\bq \cos\theta 
=\frac{(r(B)+r(C))^2 + (\rho r(B))^2 - r(C)^2}{2(r(B)+r(C))\rho r(B)} 
=\frac{r(B)+\rho^2 r(B)+ 2 r(C)}{2 \rho(r(B)+r(C))}. \label{dumpty}\eq
Also,
\bq \area(C\cap S(B))= \frac{1-\cos\theta}2\area(S(B)). \label{sat}\eq
Combining equations~(\ref{humpty}),~(\ref{dumpty}) and~(\ref{sat}),
\bq a(B,C) = \frac12-
\frac{r(B) + \rho^2 r(B)+ 2 r(C)}{4\rho(r(B)+r(C))}.\label{d2}\eq
Switching $B$ and $C$ and adding,
\bq a(B,C)+a(C,B) = 1-\frac{3 + \rho^2}{4\rho}.  \label{d3} \eq
Isn't it remarkable that $a(B,C) + a(C,B)$ does not depend on
$r(B)$ and $r(C)$?
We now choose $\rho=\sqrt 3$ to maximize the right side of equation
(\ref{d3}).  Then $a(B,C)+a(C,B) =1-\frac{\sqrt 3}2$, under the assumption
that $S(B) \cap C$ and $S(C) \cap B$ are non-empty. If $S(B) \cap C =
\emptyset$, $a(B,C) = 0$, which is greater than the negative value at the
right side of equation (\ref{d2}). As a result, $a(B,C) + a(C,B) \ge
1-\frac{\sqrt 3}2$ in the general case.  Applying this inequality to
inequality~(\ref{d1}) yields $|P| \ge |E|\left(1-\frac{\sqrt 3}2\right)$, which
gives
$$k(P) = 2|E|/|P|\le 8+4\sqrt3.$$
In conclusion, $k \le 8+4\sqrt{3}$.  By Remark~\ref{ksup}, $k(P) < k$,
establishing Theorem~\ref{ub}. \qed
\end{proof}

\begin{remark} \em In fact, $k< 8+4\sqrt 3$.  Let $B\in P$. Since each ball
$C\in E(B)$ that intersects $S(B)$ must have $r(C)\ge (\rho-1)r(B)/2$, there
is a finite bound for the number of balls $C\in E(B)$ such that $a(B,C)>0$. 
Therefore there is some $\alpha<1$ (depending on $\rho$ but not $P$) such
that
$$ \sum_{C\in E(B)} a(B,C) \le\alpha. $$
Using this inequality in place of inequality~(\ref{start}) in the above proof
would multiply the upper bound by a factor of $\alpha$.  A good estimate for
$\alpha$ would consequently strengthen Theorem~\ref{ub}.
\end{remark}

\section{The lower bound}
\begin{theorem} There exists a sequence of finite packings $\{P_n\}$ with
$$\lim_{n \to \infty} k(P_n) = 666/53.$$ \label{lb}
\end{theorem}
Observe that all questions about nerves of ball packings and average kissing
numbers are invariant under sphere-preserving transformations such as
stereographic projection from the 3-sphere $S^3$ to $\R^3$ and inversion in a
sphere.

There exists a packing $D$ in $S^3$ of 120 congruent spherical balls such that
each ball kisses exactly 12 others \cite{Coxeter}, or 720 kissing points in
total.   The existence of $D$ already implies that $k(P) > 12$ for some
packing $P$, because by Remark~\ref{ksup}, $k > k(D) = 12$.  

The proof of Theorem~\ref{lb} is a refinement of this construction.

\begin{proof}
We give an explicit description of $D$. Let $S^3$ be the unit 3-sphere in
$\R^4$ and let $\tau = {1+\sqrt{5} \over 2}$ be the golden ratio.  Choose the
centers of the balls of $D$ to be the points in the orbits of
$\frac12(\tau,1,1/\tau,0)$, $\frac12(1,1,1,1)$, and $(1,0,0,0)$ under change
of sign of any coordinate  and even permutations of coordinates.  The radius
of each ball is $18^\circ$.  We will need the following four properties of $D$,
which can be verified using the explicit description or by other means: The
12 balls that kiss a given ball have an icosahedral arrangement with 30
mutual kissing points, the centers of two kissing balls of $D$ are $36^\circ$
apart, the centers of two next-nearest balls of $D$ are $60^\circ$ apart, and
$D$ is self-antipodal.  (If $X$ is a point, set of points, or set of set of
points in $S^3$, the antipode of $X$ is given by negating all coordinates in
$\R^4$ and is denoted $-X$.)

Let $B_0 \in D$ be a ball with center $b$ and let $P_0 = D \setminus
\{B_0,-B_0\}$. The packing $P_0$ has $720-24 = 696$ kissing points and $118$
balls.  Let $R$  be the set of $12$ balls in $D$ that kiss $B_0$, and let $S$
be the unique sphere  centered at $b$ which contains the 30 kissing points
between the balls in $R$. Let $I_S:S^3 \to S^3$ be inversion in the sphere
$S$.  Observe that $S$ meets the boundary of each $B\in R$ orthogonally in a
circle (because, by symmetry, it is orthogonal to the boundary at each
kissing point), and therefore each $B\in R$ is invariant under $I_S$.  Let
$\sigma:S^3\mapsto S^3$ be the map
$\sigma(p)= I_S(-p)$.  This map $\sigma$
contracts $S^3 \setminus \{-b\}$ towards $b$, sends $-S$ to $S$, and
preserves spheres.  Because $I_S$ leaves each $B\in R$ invariant, $\sigma$
sends $-R$ to $R$.  For each $n > 0$, let
$$P_n = P_{n-1} \cup \sigma^n(P_0).$$

We claim that the sphere $S$ does not intersect any ball in $P_0 \setminus R$.
Assuming this claim, the packing $Q = P_0 \setminus (R\cup -R)$
lies between $-S$ and $S$, and $\sigma^n(Q)$ is separated from
$\sigma^{n+1}(Q)$ by $\sigma^n(S)$. Therefore each $P_n$ consists of an
alternation of layers 
$$-R, Q, \sigma(-R) = R, \sigma(Q), \sigma^2(-R),
\sigma^2(Q), \ldots, \sigma^n(-R)$$ 
such that each layer only intersects the two neighboring layers and intersects
only in kissing points.  In particular, each $P_n$ is a packing.  Moreover,
$P_{n+1}$ has $118-12 = 106$ more balls and $696-30 = 666$ more kissing
points than $P_n$ does.  Therefore
$$\lim_{n\to\infty} k(P_n) = 2 {666 \over 106} = {666 \over 53}.$$

It remains to check the claim.   Let $B_1,B_2$ be two kissing balls in $R$.
Let $b_1$ and $b_2$ be their centers and let $p$ be their kissing point.
Evidently the angular radius of $S$ is $\angle b0p$.  Using the inclusion
$S^3 \subset \R^4$ and the notation of vector calculus,
$$b_1 \cdot b_2 = b \cdot b_1 = b \cdot b_2 =  \tau/2,$$
$$b \cdot b = b_1 \cdot b_1 = b_2 \cdot b_2 = 1,$$
$$p = {b_1 + b_2\over |b_1+b_2|},$$
$$\angle b0p = \cos^{-1}\left({b \cdot (b_1 + b_2) \over |b_1+b_2|}\right)
= \cos^{-1}\left(\sqrt{2+\tau \over 5}\right) \approx 31.717^\circ.$$
On the other hand, the center of a ball in $P_0$ which is not in $R$ is at
least $60^\circ$ away from $b$, and therefore the closest point of any such
ball is at least $42^\circ$ away from $b$.  Thus, $S$ does not intersect
any such ball.
\qed \end{proof}

\bibliographystyle{plain}

\begin{thebibliography}{1}

\bibitem{C+S:book}
J.~H. Conway and N.~J.~A. Sloane.
\newblock {\em Sphere Packings, Lattices, and Groups}.
\newblock Springer-Verlag, New York, 1988.

\bibitem{Coxeter}
H.~S.~M. Coxeter.
\newblock {\em Regular Polytopes}.
\newblock Methuen, London, 1948.

\bibitem{Koebe:circpac}
P.~Koebe.
\newblock Kontaktprobleme der konformen Abbildung.
\newblock {\em Berichte Verhande.\ S\"achs.\ Akad.\ Wiss.\ Leipzig,
  Math.-Phys.\ Klasse}, 88:141--164, 1936.

\bibitem{Leech:13spheres}
J.~Leech.
\newblock The problem of the thirteen spheres.
\newblock {\em Mathematical Gazette}, 40:22--23, 1956.

\bibitem{R+S}
B.~Rodin and D.~Sullivan.
\newblock The convergence of circle packings to the Riemann mapping.
\newblock {\em J. of Differential Geometry}, 26:349--360, 1987.

\bibitem{S+W}
K.~Sch\"utter and B.~L. van~der Waerden.
\newblock Das Problem der dreizehn Kugeln.
\newblock {\em Math. Annalen}, 125:325--334, 1953.

\bibitem{Thurston:book}
W.~P. Thurston.
\newblock {\em The Geometry and Topology of 3-manifolds}.
\newblock Princeton University Notes, Princeton, New Jersey, 1982.

\bibitem{CdV:conv}
Y.~Colin~De Verdi\'ere.
\newblock Un principe variationnel pour les empilements de cercles.
\newblock {\em Invent. math.}, 104:644--669, 1991.

\end{thebibliography}

\begin{flushleft}
{\sc Department of Mathematics, University of Chicago, Chicago, IL 60637} \\
{\it E-mail address:} greg@math.uchicago.edu \\
\mbox{ } \\
{\sc Weizmann Institute of Science, Rehovot 76100, Israel} \\
{\it E-mail address:} schramm@wisdom.weizmann.ac.il
\end{flushleft}

\end{document}